\input amstex
\documentstyle{amsppt}
\magnification=\magstep1
\pageheight{9truein}
\pagewidth{6.5truein}

\NoRunningHeads
\NoBlackBoxes
\TagsOnRight
\nologo

\def \tw{\widetilde\wedge}
\def \su{\subseteq}
\def \I{(I\otimes H)}
\def \om{\omega}
\def \th{\theta}
\def \i{i_H}
\def \p{\psi^{-1}}
\def \m{(M; \Bbb R)}
\def \d{d_{HB}}

\def \o{\otimes}
\def \h{H_G^*}
\def \a{\alpha}
\def \b{\beta}
\def \hb{\widehat\beta}

\def \g{\gamma}
\def \ha{\widehat\alpha}

\def \pha{(\widehat\alpha)}

\def \ov{\overline}
\def \r{R_G\otimes}
\def \O{\Omega_{inv}}

\def \ra{\rightarrow}
\def \lra{\longrightarrow}

\topmatter
\title
The Minimal Hirsch--Brown Model via Classical Hodge Theory
\endtitle
\author
Christopher Allday and Volker Puppe
\endauthor
\abstract
In our book on cohomological methods in transformation groups the minimal 
Hirsch--Brown model was used to good effect.  The construction of the model 
there, however, was rather abstract.  Here, for smooth compact connected Lie 
group actions on smooth closed manifolds, we give a much more explicit 
construction of the minimal Hirsch--Brown model using operators from classical 
Hodge theory and the small Cartan model.
\endabstract
\endtopmatter

\document
\baselineskip=2\baselineskip

\subhead
1. Introduction
\endsubhead

In [A, P], the minimal Hirsch--Brown model is described in detail and used to
good effect.  (See, for example, [A, P] sections 1.3, 1.4, 4.4 and 4.6.)  The 
construction of the minimal Hirsch--Brown model there, however, is rather 
abstract.  Our purpose here is to give a more explicit construction of the model
 for smooth compact connected Lie group actions on closed smooth manifolds using
 operators from classical Hodge theory.  Two of our main results, Theorem (3.8)
 and Corollary (3.9), are particularly nice in view of their relation to [A2].  

In Section 2 we introduce our notation, and we give a brief description of the 
small (Cartan) model that we use for non--abelian Lie groups.  Section 3 gives 
the explicit construction of the minimal Hirsch--Brown model.  Section 4, as an
example, discusses the familiar product structure in the equivariant cohomology
 of a Hamiltonian circle action on $\Bbb C P^n$, the idea being to compute the 
 deformation of the product (as one goes from ordinary to equivariant cohomology)
 in terms of the moment map in two different ways.

\subhead
2. Notation
\endsubhead

Let $G$ be a compact connected Lie group acting smoothly on a closed smooth 
manifold $M$.  Suppose that $M$ has an invariant Riemannian metric $r$.  If one
does Hodge theory with respect to $r$, then all the usual operators, for 
example, $*$, $d^*$, $H$ (the projection onto the harmonic forms) and $G$ 
(Green's operator), are invariant.  Since the Lie group $G$ is connected, it 
acts trivially on the cohomology of $M$; and so all harmonic forms are 
invariant.  We shall let $\Omega_{inv}(M)$ or $\Omega(M)^G$ denote the cochain 
complex of invariant forms.  Thus the usual Hodge decomposition $(\alpha = 
H(\alpha)+dd^*G(\alpha)+d^*dG(\alpha))$ restricts to $\Omega_{inv}(M)$ without 
alteration.  (For an introduction to Hodge theory, see [W].)  

In this paper we shall always assume that the Lie group, $G$, the manifold, $M$,
 and the Riemannian metric are as in the paragraph above.

When $G=T^k$, the $k$--torus, we shall use the Cartan model to compute the 
equivariant cohomology, $H^*_G(M; \Bbb R)$.  Let $R_G = H^*(BG; \Bbb R) = 
\Bbb R [t_1, \dots, t_k]$, where each $t_i$ has degree $2$.  There is a 
canonical association between $t_1, \dots, t_r$ and a basis $c_1, \dots, c_r$ of
 $\frak g$, the Lie algebra of $G$.  And, to each $c_j$, there is associated a 
vector field $V_j$ on $M$, via the exponential map, $\frak g \rightarrow G$, and
the group action.  The Cartan model is then $R_G \otimes \Omega_{inv}(M)$ with 
differential $d_G = I\otimes d - \partial$, where $\partial = \sum_{j=1}^kt_j
\otimes i_j$ and $i_j = i_{V_j}$, the inner product, for $1\le j\le k$.  We 
usually abbreviate $I\otimes d$ as $d$; and so $d_G = d-\partial$.  (See [A, B]
or [G, S] for details.)   

For non--abelian $G$, it is convenient here to use the small Cartan model.  
This was proposed by a number of authors (see, for example, [G, K, M]), and 
it has been shown to be equivalent to the Cartan model by Alekseev and 
Meinrenken ([A, M]).  Formally, the small model looks just like the above, 
that is, $R_G\otimes \Omega_{inv}(M)$ with differential $d-\partial$, where 
$R_G = H^*(BG; \Bbb R) = \Bbb R[t_1, \dots, t_k]$ and $\partial = 
\sum\limits_{j=1}^kt_j\otimes i_j$.  So $k$ is the rank of $G$.  Each $i_j$ 
is now, 
however, the inner product with a multivector corresponding to an element 
$c_j \in (\wedge\frak g)^G$, where $\{c_1, \dots, c_k\}$ is a basis of 
$(\wedge\frak g)^G$.  Each $c_j$ has positive odd degree, but the degrees 
may be more than one.  The polynomial generator $t_j$ corresponds to $c_j$ 
via a canonical transgression.  So $\deg$ $t_j = \deg$ $c_j + 1$; and, as an 
operator, $\deg$ $i_j = 1 - \deg$ $t_j$.  For example, when $G = SU(2)$, 
$k=1$, $\deg$ $t_1 = 4$ and $\deg$ $i_1 = -3$.  For details of the 
construction of the small model and its equivalence to the Cartan model, 
see [A, M].

Because, typically, some or all $i_j$ in the small model are inner products 
with multivectors, $d_G = d-\partial$ is not a derivation in the non--abelian 
case:  so the product in the small model is not the obvious one.  Again, 
see [A, M] for more details.  This, however, is not a problem here since 
the product in the minimal Hirsch--Brown model is not the obvious one either.

There are two operators which play an important role in our description of 
the minimal Hirsch--Brown model and its relation to the small model (which 
coincides with the Cartan model when $G$ is a torus).

\demo{Definition (2.1)}  Let $P = (I\otimes d^*G)\partial$ and $Q = \partial
(I\otimes d^*G)$, where $G$ is Green's operator.  More briefly, $P = d^*G
\partial$ and $Q = \partial d^*G$.
\enddemo

Since $d^*$ and $G$ commute, $PQ = QP = 0$.

\subhead 3.  The Minimal Hirsch--Brown Model  \endsubhead

We begin with the following abbreviations.

\demo {Definitions (3.1)}
(1) On $R_G\otimes \Omega_{inv}(M)$, let $\varphi = I-P$ 
and $\psi = I-Q$, where $P = d^*G\partial$ and $Q = \partial d^*G$, as above.

(2)  In $\O(M)$, let $\Cal H$ be the subspace of harmonic forms, $B=im(d)$ be the boundaries and $E=im(d^*)$ be the coboundaries.  So, in the Hodge
decomposition, $\O(M) = \Cal H \oplus B \oplus E$.
\enddemo

Restricting $\psi$ to $R_G\otimes B$, gives the following commutative diagram.

\demo {Lemma (3.2)} 
The following diagram commutes, where the top arrow is the inclusion.

$$
\CD
\r E      @>>>     \r \Omega_{inv}(M)   \\
@VdVV                  @VVd_GV \\
\r B      @>>\psi>    \r \Omega_{inv}(M)
\endCD 
$$
\enddemo

\demo {Proof} 
$\psi d|R_G\otimes E = (d-\partial d^*dG)|R_G\otimes E = d_G|R_G\otimes E$, 
because $d^*dG$ is the identity on $E$.
\enddemo

(To put it another way, $d_Gd^*G|R_G\otimes B = \psi|R_G\otimes B.)$  \qed

Now we define a new differential, $\overline D$ on $R_G\otimes\Omega_{inv}$ and
the Hirsch--Brown differential, $d_{HB}$, on $R_G\otimes \Cal H$.                                                  

\demo{Definitions (3.3)}
(1)  On $R_G\otimes\Omega_{inv}(M)$, set $\overline D = \psi^{-1}d_G\psi$.  

(2)  On $\r \Cal H$,
put $d_{HB} = (I\o H)\ov D|R_G\o \Cal H$.
\enddemo

It is clear that $\ov D^2 = 0$; that $d_{HB}^2 = 0$, too, follows from the next
lemma.

\demo {Lemma (3.4)}
The following diagram commutes.
$$
\CD
\r \O (M)  @>\psi^{-1}>>   \r\O (M)   @>I\o H>>     \r \Cal H      \\
@Vd_GVV    @V\ov DVV       @VVd_{HB}V                               \\
\r\O(M)    @>>\psi^{-1}>   \r\O(M)    @>>I\o H>     \r \Cal H       
\endCD
$$
\enddemo
 
\demo {Proof}
It is enough to show that, for any $\alpha\in B\oplus E$, $(I\o H)\ov D(\alpha) = 
0$.  This follows from the next lemma.  \qed
\enddemo

\demo {Lemma (3.5)}  
(1)  On $\r\Cal H$, $\ov D = -\psi^{-1}\partial = -\partial\varphi^{-1}$.  

(2)  On $\r B$, $\ov D = 0$.

(3)  On $\r E$, $\ov D = d$.
\enddemo

\demo {Proof} 
(1)  On $\r \Cal H$, $d = 0$ and $Q = 0$  :  so $\ov D = \psi^{-1}d_G\psi = 
\psi^{-1}d_G = -\psi^{-1}\partial$.  And $Q\partial = \partial P$.

(2)  From the proof of Lemma (3.2), $\psi|\r B = d_Gd^*G|\r B$.
Hence $d_G\psi|\r B = 0$.

(3)  On $\r E$, $Q = 0$; and so, on $\r E$,
$\ov D = \psi^{-1}d_G = \psi^{-1}(\psi d)$, by Lemma (3.2)  \qed
\enddemo

\demo {Definition (3.6)}  
The differential $R_G$--module, $(\r \Cal H, d_{HB}) = (H^*(BG; \Bbb R)\otimes 
H^*(M; \Bbb R)$, $d_{HB}$) is called the minimal Hirsch--Brown model for $H^*_G
(M; \Bbb R)$.  That it computes $H^*_G(M; \Bbb R)$ follows from the next lemma.
\enddemo

\demo {Lemma (3.7)} 
$H(\r \Cal H, d_{HB}) \cong H^*_G(M; \Bbb R)$, where the $H$ on the left means 
(co)homology with respect to the differential $d_{HB}$.

Indeed, $(I\o H)\psi^{-1} = (1\o H)(I-Q)^{-1}$ is a homotopy equivalence of 
differential $R_G$--modules.
\enddemo

\demo {Proof}  
Since $\psi$ is an isomorphism, $\psi^{-1}$ induces an isomorphism on cohomology
.  And, by Lemma (3.5) (2) and (3), $\r (B\oplus E) = \ker(I\o H)$ is acyclic with 
respect to $\ov D$.  
So $I\o H$ also induces an isomorphism in cohomology.  This proves the first 
statement of the lemma.
\enddemo

Since all the differential $R_G$--modules involved are free, so is the mapping 
cone of $(I\o H)\psi^{-1}$.  Thus the second statement of the lemma follows from
 [A, P], Remark (B.1.10), Proposition (B.1.11) and Proposition (B.1.7)  \qed

The Hirsch--Brown differential can be written in a very useful way as we show 
next.

\proclaim {Theorem (3.8)} 
On $\r\Cal H$, 
$$
d_{HB} = (I-P)d_G(I-P)^{-1} = \varphi d_G\varphi^{-1}.
$$
\endproclaim

\demo {Proof} 
Let $a\in \r \Cal H$.  By Lemma (3.5)(1), 
$$ 
\split
d_{HB}(a) 
&= (I\otimes H)\ov D(a) = -(I\o H)\partial \varphi^{-1}(a) \\
&= -\partial\varphi^{-1}(a)+(I\o\Delta G)\partial\varphi^{-1}(a),
\endsplit
$$ 
by the Hodge Decomposition Theorem.
\enddemo
          
In general, however, 
$$
(I\o \Delta G)\partial = d^*Gd\partial +dd^*G\partial = dP-Pd,
$$
where, as usual, we have abbreviated $(I\o d^*dG)\partial, (I\o d)P$, et 
cetera, by $d^*dG\partial$, $dP$, \break et cetera.  I.e., $[d, P] = \Delta 
G\partial$.
$$
\split
\text{So\ } \Delta G\partial\varphi^{-1}(a) 
&= dP\varphi^{-1}(a)-Pd\varphi^{-1}(a)\\
&= d(I-\varphi)\varphi^{-1}(a)-Pd\varphi^{-1}(a)\\
&= d\varphi^{-1}(a)-Pd\varphi^{-1}(a), \text{\ since\ } d(a) = 0.
\endsplit
$$
Finally, then, $d_{HB}(a) = -\partial\varphi^{-1}(a)+\varphi d\varphi^{-1}(a) = 
\varphi d_G\varphi^{-1}(a)$, since $P\partial = 0$.  \qed

\demo {Corollary (3.9)}  
The following diagram commutes, where $i_H$ is the inclusion.
$$
\CD
\r \Cal H      @>i_H>>       \r\O (M)     @>\varphi^{-1}>>     \r\O(M)  \\
@Vd_{HB}VV     @VV\varphi d_G\varphi^{-1}V      @VVd_GV                    \\
\r\Cal H       @>>i_H>       \r\O(M)      @>>\varphi^{-1}>     \r\O(M) 
\endCD
$$
\enddemo
Furthermore, $(I \otimes H)\psi^{-1}\varphi^{-1}i_H = I$.

\demo {Proof}
Since $PQ = QP = 0$, $\p\varphi^{-1} = I+P+Q+(P+Q)^2+\cdots$.  And $HP = 0$ 
and $QH = 0$.  So $\I\p\varphi^{-1}\i = I.$  \qed
\enddemo

\demo {Remark (3.10)}
Corollary (3.9) shows that $\I\p$ is a fibration and $\varphi^{-1}\i$ is a 
cofibration by [A, P], Proposition (B.1.5.) and Proposition (B.1.4.).
\enddemo

Theorem (3.8) also reproves the main result of [A2].  Let $i:  M\ra M_G$ be 
the inclusion of a fibre in the Borel construction bundle $M_G\lra BG$.

\demo {Corollary (3.11)}
Suppose that $i^*:  H^*_G\m\ra H^*\m$ is surjective.  Let $\alpha\in\O(M)$ be a
harmonic form.  Then $\alpha$ has a canonical equivariant extension; namely
$$
d_G(I-P)^{-1}(\alpha) = 0\,.
$$
\enddemo

\demo {Proof}
Since $i^*$ is surjective, it follows that 
$$
H^*_G\m\cong\r H^*\m
$$
as a $R_G$--module.  Hence $\d = 0$.  Thus, by Theorem (3.8), $d_G(I-P)^{-1}
(\alpha) = 0$, for all $\alpha\in\Cal H$.  \qed
\enddemo

The operators used above, $P, Q, H, \Delta, G$, for example, are not 
multiplicative, and nor is the small model itself when the Lie group is 
non--abelian.  So, in general, even for torus actions, it is not easy to 
describe the product structure in the minimal Hirsch--Brown model; and, even 
when $\d$ is zero, the product on $\r \Cal H$ is usually twisted.  Nevertheless,
 for torus actions when $\d$ is zero, we can describe the product in $H^*_G\m$.

\demo {Definition (3.12)}  
(1)  As in [A2], we shall use $CEF$ to mean that there is a cohomology 
extension of the fibre, that is, 
$$
i^*:  \h\m\lra H^*\m
$$
is surjective.  When $G$ is a torus, this implies that 
$$
j^*:  \h\m\lra \h(M^G; \Bbb R)
$$
is injective where $j:  M^G\ra M$ is the inclusion of the fixed point set.  
(See, e.g., [A, P], Section (3.1).)  And, for any compact connected $G$, as 
noted above, $CEF$ implies that
$$
\h\m\cong H^*(BG; \Bbb R)\o H^*\m
$$
as $R_G$--modules.  Thus $CEF$ implies that $\d = 0$. 

Note that, when using cohomology with coefficients in an abelian group that is
not a field, then one says that there is a $CEF$ if $i^*$ has a right inverse.
See, e.g., [S], Chap.5, Sec.7. This is, of course, equivalent to the surjectivity
of $i^*$ when using field coefficients. 

(2)  Let $\tw$ denote the product in the minimal Hirsch--Brown model.  
In particular, for $\a, \b\in\Cal H\su\O(M)$, $\a\tw\b$ is the product of $\a$ 
and $\b$ in $\r \Cal H$, whereas, of course, $\a\wedge\b$ is the product in 
$\O(M)$.    
 
(3)  For $\a\in\O(M)$, abbreviate $(I-P)^{-1}(\a) = \varphi^{-1}(\a)$ by $\ha$.
\enddemo

We now have the following description of the cup--product in $\h\m$ in the
$C
EF$ case when $G$ is a torus.  Of course, 
$$
\h\m\cong\r\Cal H;
$$
and so it is enough to describe $\tw$.  Indeed, since $H^*_G(M; \Bbb R)$ is a 
$R_G$--algebra, it is enough to describe $\tw$ on $\Cal H$.

\demo {Proposition (3.13)}
Suppose that $G$ is a torus and that there is a $CEF$.  Then, for $\a, 
\b\in\Cal H$,
$$
\a\tw\b = \I(1-Q)^{-1}\big(\ha\hb\,\big).
$$
\enddemo

\demo {Proof}
Let $\th = \I(1-Q)^{-1}$.  Since $\d = 0$, $\ha$ and $\hb$ are cycles in 
$(\r\O(M), d_G)$ by Corollary (3.11).  And, since $G$ is a torus, $(\r\O(M), 
d_G)$ is the Cartan model and, thus, multiplicative.  Hence the product $\ha\hb$
 in $\r\O(M)$ represents $[\ha]\big[\hb\,\big]$, the product in $\h\m$.  And, 
although $\th$ is not multiplicative, $\th^*$ is an isomorphism; and so we have 
the following.
$$
\split
\a\tw\b 
&= \th^*([\ha])\th^*\big(\big[\hb\,\big]\big), \text{since\ } \d = 0 \\
&= \th^*\big([\ha]\big[\hb\,\big]\big) = \th^*\big(\big[\ha\hb\,\big]\big)\\
&= \big[\th\big(\ha\hb\,\big)\big] = \th\big(\ha\hb\,\big), \text {again 
because\ } \d = 0 \qed
\endsplit
$$
\enddemo
\demo {Remarks (3.14)}
(1)  Since $\ha = (I-P)^{-1}(\a)$, for $\a, \b\in\Cal H$, under the conditions 
of Proposition (3.13), 
$$
\a\tw\b = H(\a\wedge\b) \text{ modulo\ } \ov R_G\o\Cal H,
$$
where $\ov R_G$ is the augmentation ideal of elements of positive degree in 
$R_G$.  And, of course, $H(\a\wedge\b)$ is the product of $\a$ and $\b$ in 
$\Cal H\cong H^*\m$.

(2)  If $M$ is a closed symplectic manifold and the action of $G$ (any compact 
connected Lie group) is symplectic, then there is a $CEF$ if and only if the 
action is Hamiltonian.  This follows largely from results of Frankel ([F]).

(3)  An argument very similar to the proof of Theorem (3.8) shows that for any 
$\a\in\O(M)$,
$$
(I-P)d_G(I-P)^{-1}(\a) = (I\o H)d_G(I-P)^{-1}(\a)+d\a.
$$
(Briefly, $Hd_G\pha = -H\partial\pha = -\partial\pha+\Delta G\partial\pha 
= -\partial\pha+dP\pha-Pd\pha = -\partial\pha+d(\ha-\a)-Pd\pha 
=d_G\ha-Pd_G\pha-d\a.)$

(4)  Under the conditions of Proposition (3.13), from Corollary (3.9), it 
follows similarly that, for any $\a, \b\in\Cal H$, there is $\g\in\r\O(M)$ such 
that 
$$
\varphi^{-1}\i(\a\tw\b) = \ha\hb + d_G\g.
$$
  
(5)  Similar results hold for products of three or more elements.
\enddemo

\subhead   
4.  An Example
\endsubhead
Let $M$ be a closed symplectic $2n$--manifold with symplectic form $\om$.  
Suppose that a compact connected Lie group, $G$, is acting on $M$ in a 
Hamiltonian way.  Then we may choose an invariant Riemannian metric on $M$ that 
is compatible with $\om$.  See, e.g., [M, S], Lemma 5.49.  So, if $r$ is the 
metric, and $V_1$ and $V_2$ are any two vector fields on $M$, then $r(V_1, V_2)
= \om(V_1, JV_2)$, where $J$ is an invariant compatible almost--complex 
structure on $M$.  It follows that $\om^j$ is harmonic for $0\le j\le n$; and 
$$
\ast\left(\frac{\omega^j}{j!}\right) = \frac{\omega^{n-j}}{(n-j)!}\,, 
$$
for $0\le j\le n$.  In particular, $\frac{\omega^n}{n!}$ is the volume form.  

As remarked above, it follows from the results of Frankel ([F]) that, in the 
Hamiltonian case, $M$ has a $CE F$; and so, in the minimal Hirsch--Brown 
model, $\d = 0$.  Thus the remaining problem is to determine the product 
structure in $\h\m$.  In this section we shall do this in the familiar situation
 where $G=S^1$ and $M = \Bbb CP^n$.  The results are not new, although they may 
be assembled in a somewhat novel way.

First, however, consider a Hamiltonian action of $G = S^1$ on any closed 
symplectic manifold $(M, \om)$.  Let $\mu$ be the moment map; and suppose that 
$\mu$ has been chosen to have average value zero on $M$:  i.e., $\int_M\mu
\frac{\om^n}{n!} = 0$.  Let $V$ be the vector field defined by the circle action
:  so, for any $x\in M$,
$$
V_x = \frac d{du}\exp(2\pi iu)x|_{u=0}.
$$
In the Cartan model, then, the differential $d_G = d-ti_V$, where $t\in H^2(BG;
\Bbb R)$ is the polynomial generator.  In the Hodge decomposition $\mu = d^*dG
\mu$, since the harmonic part, $H(\mu)$, is the average value.  Thus $P(\om) = 
t\mu$, because $d^*Gi_V(\om) = d^*Gd\mu = \mu$.  Hence $\widehat\om = \om+t\mu$,
 the standard equivariant extension of $\om$.

Now let $M = \Bbb CP^n$ with symplectic form $\om$ and Hamiltonian action of 
$G = S^1$.  Let $\mu$ be the moment map; but we shall not assume that $H(\mu) =
0$.  Let $\overline w = [\om+t\mu]_G\in H^2_G\m$.  The product structure in
$\h\m$ is completely determined by expressing $\overline w^{n+1}$ in terms of 
lower powers of $\overline w$.  Let $\overline w^{n+1} = \sum\limits_{i=1}^
{n+1}c_i\overline w^{n+1-i}t^i$, where each $c_i \in\Bbb R$.  One way to find 
the $c_i$ is the following.

For $j\ge 0$, $\overline w^{n+1+j} = \sum\limits_{i =1}^{n+1}c_i\overline w^
{n+1+j-i}t^i$.  So integrating over the fibre, $M$, in the Borel construction 
bundle $M_G\ra BG$, gives 
$$
\left ( \matrix n+1+j\\1+j\endmatrix\right)t^{1+j}\int_M\mu^{1+j}\om^n = 
\sum_{i=1}^{1+j}c_i
\left ( \matrix n+1+j-i\\1+j-i\endmatrix\right)t^{1+j}\int_M\mu^{1+j-i}
\om^n\,.  
$$
$$
\text {So\ } \pmatrix n+1+j\\1+j\endpmatrix H(\mu^{1+j}) = \sum^{1+j}_{i=1}c_i 
\pmatrix n+1+j-i\\1+j-i\endpmatrix H(\mu^{1+j-i}).
$$
Since this holds for all $j\ge 0$, one can easily solve for each $c_i$ in terms
of the average values of the powers of $\mu$.  For example, putting $j=0$, $c_1
 = (n+1)H(\mu)$; and, putting $j=1$, 
$$
c_2 = \left ( \matrix n+2\\2\endmatrix\right)H(\mu^2)-(n+1)^2H(\mu)^2.
$$
Equally, one can solve for each $H(\mu^j)$ in terms of $c_1, \dots, c_j$.  
This is reasonable because there are other familiar ways to find the $c_is$.  
Let the fixed point set $M^G = \bigcup\limits^s_{i=1}F_i$, where the component 
$F_i$ 
has dimension $2r_i$.  By the equality of Euler characteristics, 
$\sum\limits_{i=1}^s(r_i+1) = n+1$.  Let $\nu_i$ be the value of $\mu$ on 
$F_i$; and let $\mu_j = \nu_i$ for $\sum\limits^{i-1}_{k=1}(r_k+1)+1\le j\le 
\sum\limits^i_{k=1}(r_k+1)$.  So the distinct values of $\mu$ appear with 
multiplicity, each $\nu_i$ appearing with multiplicity $r_i+1$.  (If $M^G$ 
is finite, then $s=n+1$, and $\mu_i = \nu_i$ for $1\le i\le n+1$.)  In terms 
of these values we have
$$
\prod^s_{i=1}(\overline w-\nu_it)^{r_i+1} = \prod^{n+1}_{i=1}(\overline w-
\mu_it) = 0.
$$
This follows from the Localization Theorem of Borel, Hsiang and Quillen.  For 
details of this example see [H], Theorem (IV.3) or, e.g., [A1], Example (3.12). 
 Thus, for $1\le i\le n+1$, $c_i = (-1)^{i+1}\sigma_i$, where $\sigma_i$ is the
 $i$th elementary symmetric polynomial in $\mu_1, \dots, \mu_{n+1}$.

Now suppose that $M^G$ is finite:  so $s = n+1$ and each $r_i = 0$.  Let 
$U_i = \prod\limits_{j\ne i}(\overline w-\mu_jt)$.  So $U_i$ restricts to 
$\prod\limits_{j\ne i}(\mu_i-\mu_j)t^n$ at $F_i$ and zero at all the other 
fixed points.  Let the equivariant Euler class at $F_i$ be $\varepsilon_it^n$, 
normalized so that $\varepsilon_i$ is an integer (the product of the weights).
  Integrating $U_i$ over the fibre gives
$$
\int_M\om^n = \frac 1{\varepsilon_i}\prod_{j\ne i}(\mu_i-\mu j)\tag4.1
$$
by the Integration Formula ([A, B], (3.8)).  (See [B], Chapter VIII, Theorem 5.5
 (based on the original example of W.--Y. Hsiang), [P] for many related results,
 or [A1], Example (4.16) for an elementary treatment.)

Meanwhile the Duistermaat--Heckman formula gives
$$
\int_Me^{\mu t}\frac{\om^n}{n!} = \sum^{n+1}_{i=1}\frac{e^{\mu_it}}
{\varepsilon_it^n}.
$$
Thus
$$
\left ( \matrix n+j\\j\endmatrix\right)H(\mu^j) = \sum^{n+1}_{i=1}
\frac{\mu_i^{n+j}}{\prod_{j\ne i}{(\mu_i-\mu_j)}}.
$$

The last formula is homogeneous in $\mu$ and, hence, not sensitive to such 
matters as how on parametrizes the circle $(\exp(2\pi it)$ versus $\exp(it))$, 
what sign convention one uses for $\mu(d\mu = i_V(\om)$ versus $d\mu = 
-i_V(\om))$ or the sign of $t(d_G = d-ti_V \text{\ versus\ } d_G = d+ti_V)$.

Given a particular linear action, one can use (4.1) to find the $\mu_is$, and, 
hence, the $c_is$.  For example, let $S^1$ act on $\Bbb CP^2$ by $z[z_0, z_1, 
z_2] = [z_0, z^az_1, z^bz_2]$, where $a$ and $b$ are integers such that $0<a<b$.
  Let $\int_M\om^2 = A$.  Choose $\mu$ so that $H(\mu) = 0$.  Then one gets 
$$
\align
6H(\mu^2) &= \frac A3(a^2-ab+b^2) = c_2\\
\text{and\ }10H(\mu^3) &= \frac{A\sqrt A}{27}\,(2a^3-3a^2b-3ab^2+2b^3) = c_3\,.
\endalign
$$
Thus
$$
\overline w^3 = \frac 13(a^2-ab+b^2)A\overline w t^2+\frac 1{27}(2a^3-3a^2b-
3ab^2
+2b^3)A\sqrt A\,t^3\,.
$$
$(c_3 = \mu_1\mu_2\mu_3 = \frac {-A\sqrt A}{27}(a+b)(2a-b)(2b-a)$, for example.)

\enddocument